\documentclass[12pt]{amsart}%
\usepackage{eurosym}
\usepackage{amssymb}
\usepackage{amsmath}
\usepackage{amsfonts}
\usepackage{graphicx}%
\newtheorem{theorem}{Theorem}[section]

\newtheoremstyle{notauto}{}{}{\itshape}{}{\bfseries}{.}{0.5em}{\thmnote{#3}}
\theoremstyle{notauto}

\theoremstyle{definition}

\theoremstyle{remark}
\newtheorem{remark}[theorem]{Remark}

\renewcommand{\geq}{\geqslant}
\renewcommand{\leq}{\leqslant}

\begin{document}
\title{A short note \\on the noncoprime regular module problem}
\author{G\"{u}l\.{I}n Ercan$^{*}$}
\address{G\"{u}l\.{I}n Ercan, Department of Mathematics, Middle East Technical
University, Ankara, Turkey}
\email{ercan@metu.edu.tr}
\author{\.{I}sma\.{I}l \c{S}. G\"{u}lo\u{g}lu}
\address{\.{I}sma\.{I}l \c{S}. G\"{u}lo\u{g}lu, Department of Mathematics,
Do\u{g}u\c{s} University, Istanbul, Turkey}
\email{iguloglu@dogus.edu.tr}
\thanks{$^{*}$Corresponding author}
\thanks{This work has been supported by the Research Project T\"{U}B\.{I}TAK 114F223.}
\maketitle

\begin{abstract}
We consider a special configuration in which a finite group $A$ acts by
automorphisms on the finite group $G$, and the semidirect product $GA$ acts on
the vector space $V$ by linear transformations; and discuss the existence of
the regular $A$-module in $V_{_{A}}$.

\end{abstract}

\label{s1}

\section{Introduction}

Let $A$ be a finite group which acts faithfully on the vector space $V$ by
linear transformations. We say \textquotedblleft$A$ has a regular orbit on
$V$" if there is a vector $v$ in $V$ such that $C_{A}(v)=1$. In this case, the
$A$-orbit containing $v$ is called a regular $A$-orbit. Furthermore, $V$
contains the regular $A$-module if a regular $A$-orbit happens to be linearly
independent. More generally if $A$ acts by linear transformations on the
vector space $V$ (not necessarily faithfully), then we say that $A$ has a
regular orbit on $V$ or $V$ contains the regular $A$-module if $A/C_{A}(V)$
does the same.

While studying the structure of a finite solvable group $G$ admitting a
certain group of automorphisms $A$, we are often forced to study $A$-invariant
chief factors $V$ of $G$ together with the action of the semidirect product
$(G/C_{G}(V))A$ on $V.$ It turns out to be rather efficient to know that $V$
contains the regular $A$-module or at least a regular $A$-orbit. Not all
groups act with regular orbits although many interesting and rich classes do,
especially under the additional assumptions of coprimeness that
$(|G|,|A|)=1=(|V|,|GA|).$ There has been extensive research about the
existence of regular orbits such as \cite{Be}, \cite{Fle}, \cite{Gow},
\cite{Har}, \cite{Turull}, \cite{Tur2} in the case of coprimeness and
\cite{Car}, \cite{Esp}, \cite{Esp1}, \cite{Ya1}, \cite{Ya2} in the noncoprime
case. All the results concerning a nilpotent $A$ are culminating in Theorem
1.1 in \cite{Ya2} which can be reformulated as follows:\newline

\textit{ Let $G$ be a finite solvable group admitting a nilpotent group $A$ as
a group of automorphisms. Suppose that $C_{O_{p}(A)}(G)=1$. Let $V$ be a
finite faithful $kGA$-module over a field $k$ of characteristic $p$ not
dividing the order of $G$. Then $A$ has at least one regular orbit on $V$ if
$A$ involves no wreath product $\mathbb{Z}_{2}\wr\mathbb{Z}_{2}$ and involves
no wreath product $\mathbb{Z}_{r}\wr\mathbb{Z}_{r}$ for $r$ a Mersenne prime
when $p=2.$}\newline

In the present paper we prove a theorem which concludes the existence of a
regular module without the coprimeness condition the prototype of which is
Theorem 1.5 in \cite{Turull}. This theorem was improved as Theorem B in
\cite{Esp1} in case where the group $GA$ is of odd order. For the convenience
of the reader, we formulate the main conclusion of Theorem 1.5 in a way
suitable to emphasize the similarities and differences between this theorem
and Theorem B in \cite{Esp1} and our result.\newline

\textit{Let $PRA$ be a finite group where $P$ is a $p$-group and $R$ is an
$r$-group for distinct primes $p$ and $r$ not dividing the order of $A$ such
that $P \lhd PRA$ and $R \lhd RA$. Assume that the following are
satisfied:}\newline\textit{(a) $P$ is an extraspecial $p$-group for some prime
$p$ where $Z(P)\leq Z(PRA)$ and $C_{A}(P)=1$;}\newline\textit{(b) $\bar
{R}=R/R_{0}$ is of class at most two and of exponent $r$ where $R_{0}%
=C_{R}(P)$. Suppose that $|C_{A}(\bar{R}/\Phi(\bar{R})|$ is either a prime or
$1$;}\newline \textit{(c) $A/C_{A}(\bar{R}/\Phi(\bar{R})$ has a regular orbit
in its action on $\bar{R}/\Phi(\bar{R})$;\newline if $C_{A}(\bar{R}/\Phi
(\bar{R})\ne1,$ $[C_{A}(\bar{R}/\Phi(\bar{R}),P]\ne P$ and $p=2$, assume that
$|C_{A}(\bar{R}/\Phi(\bar{R})|$ is not a Fermat prime.}\newline\textit{Let
$\chi$ be a complex $PRA$-character such that $\chi_{_{P}}$ is faithful. Then
$\chi_{_{A}}$ contains the regular $A$-character.}\newline

Namely we obtain the following:\newline

\textbf{Theorem} \textit{Let $PRA$ be a finite group where $P$ is a $p$-group
and $R$ is an $r$-group for distinct primes $p$ and $r$ such that $P \lhd PRA$
and $R \lhd RA$. Assume that the following are satisfied:\newline(a) $P$ is an
extraspecial $p$-group for some prime $p$ where $Z(P)\leq Z(PRA)$ and
$C_{A}(P)=1$;\newline(b) $R/R_{0}$ is of class at most two and of exponent
dividing $r$ where $R_{0}=C_{R}(P)$ and $A_{0}=C_{A}(R/R_{0})= 1;$\newline(c)
$A=A_{p} \times A_{r} \times A_{\{p,r\}^{\prime}}$ where its Sylow
$r$-subgroup $A_{r}$ and Sylow $p$-subgroup $A_{p}$ are both cyclic and
$A_{\{p,r\}^{\prime}}$ acts with regular orbits on $R/\Phi(R)$,\newline(d) if
$p=2$ then $r$ is not a Fermat prime.\newline Let $\chi$ be a complex
$PRA$-character such that $\chi_{_{P}}$ is faithful. Then $\chi_{_{A}}$
contains the regular $A$-character.}\newline

Notice that both $p$ and $r$ are allowed to divide the order of $A$.\newline

All groups considered in this paper are finite and the notation and
terminology are standard.

\section{Existence of regular orbits}

In this section we present a result due to Dade \cite{Dade} on the existence
of regular orbits which will be applied in the proof of our theorem. \newline

\textbf{Proposition} \textit{Let V be a faithful $kA$-module over a finite
field $k$ of characteristic $p$. Assume that $A=B\times C$ where $B$ is a
cyclic $p$-group and $C$ is a $p^{\prime}$-group which has a regular orbit on
every $C$-invariant irreducible section of $V$. Then $A$ has a regular orbit
on $V$. }

\begin{proof}
Let $V_{C}=W_{1}\oplus\cdots\oplus W_{\ell}$ be the decomposition of $V$ into
its $C$-homogenous components. As $B$ and $C$ commute, each $W_{i}$ is
$A$-invariant. Therefore it suffices to prove that $A$ has a regular orbit on
$W_{i}$, for each $i=1,\ldots,\ell$. To see this let $w_{i}\in W_{i}$ be such
that $C_{A}(w_{i})=C_{A}(W_{i})$ for $i=1,\ldots,\ell$. If $v=w_{1}%
+\cdots+w_{\ell}$, then
\[
C_{A}(v)=\bigcap_{i=1}^{k}C_{A}(w_{i})=\bigcap_{i=1}^{k}C_{A}(W_{i}%
)=C_{A}(V)=1.
\]
Thus we may assume that $\ell=1$, that is, $V_{C}$ is homogeneous. Let $X$ be
the irreducible $kC$-module which appears in $V_{C}$ and let $B=\left\langle
\alpha\right\rangle $. Then we have $kB=k[\alpha-1]$. Set $R_{j}%
=kB/\left\langle (\alpha-1)^{j}\right\rangle $, for $j=1,\ldots,p^{n}$ where
$p^{n}=|\alpha|$. Note that $R_{j}$ is an indecomposable $kB$-module of
dimension $j$ for each $j$ and these are the only indecomposable $kB$-modules
by Theorem VII.5.3 in \cite{HB}. Then the decomposition of the $kA$-module $V$
into indecomposable $kA$-modules can be given as
\[
V\cong(X\otimes R_{j_{1}})\oplus\cdots\oplus(X\otimes R_{j_{m}})\cong
X\otimes(\bigoplus_{i=1}^{m}R_{j_{i}})
\]
for some $j_{1},\ldots,j_{m}$ in $\{1,\ldots,p^{n}\}$. To simplify the
notation we set $U=\bigoplus_{i=1}^{m}R_{j_{i}}$. The group $C$ has a regular
orbit on $X$ by hypothesis, that is, there is $x\in X$ such that
$C_{C}(x)=C_{C}(X)=1$. We shall observe that $B$ has a regular orbit on $U$ :
As a consequence of the faithful action of $A$ on $V$, $B$ acts faithfully on
$U$. Hence there is at least one indecomposable component, say $R_{j_{i}}$, on
which $B$ acts faithfully, since $B$ is cyclic. Let
\[
R_{j_{i}}=U_{1}\supset U_{2}\supset\cdots\supset U_{s}=0
\]
be a $B$-composition series of $R_{j_{i}}=U_{1}$. Each factor $U_{i}/U_{i+1}$,
$i=1,\ldots,s-1$, is isomorphic to the trivial module of dimension $1$. Hence
$s-1=dimU_{1}=j_{1}$ and $[U_{1},\underbrace{\alpha,\ldots,\alpha}%
_{j_{1}-times}]=0$. It follows that $dimU_{1}=j_{1}\geq p^{n-1}+1$, because
otherwise $(\alpha-1)^{p^{n-1}}=0$ on $U_{1}$ and hence $\alpha^{p^{n-1}}$ is
trivial on $U_{1}$, a contradiction. Pick an element $u$ from $U_{1}-U_{2}$.
If $C_{B}(u)\neq1$, then $\alpha^{p^{n-1}}$ acts trivially on $u$, whence the
degree $j_{1}$ of the minimum polynomial of $\alpha$ on $U_{1}$ is at most
$p^{n-1}$ . But then $p^{n-1}+1\leq j_{1}\leq p^{n-1}$, which is impossible.
This yields that $C_{B}(u)=1=C_{B}(U)$. As a consequence, $B$ has a regular
orbit on $U$. We are now ready to complete the proof of the theorem. Let $a\in
C_{A}(x\otimes u)$. Then $a=b+c$ for some $b\in B$ and $c\in C$. As
$c\in\left\langle a\right\rangle $, we have $(x\otimes u)c=xc\otimes
u=x\otimes u$ and hence $xc=x$. That is, $c\in C_{C }(x)=C_{C}(X)$. Similarly,
we observe that $b\in C_{B}(u)=C_{A}(U)$. Consequently, we have $a\in
C_{A}(X\otimes U)$ and hence the equality $C_{A}(x\otimes u)=C_{A}(X\otimes
U)$ holds. It follows that $A$ has regular orbit on $V$, as claimed.
\end{proof}

\begin{remark}
The above proposition cannot be extended to abelian $O_{p}(A)$ as the
following example shows: Let $V$ be an elementary abelian group of order
$p^{3}$ with a basis $\{v_{1},v_{2},v_{3}\}$ and $A$ an elementary abelian
group of order $p^{2}$ of automorphisms of $V$ generated by $\{a_{1},a_{2}\}$
with the action $v_{1}^{a_{1}}=v_{1}^{a_{2}}=v_{1},v_{2}^{a_{1}}=v_{1}%
v_{2},v_{2}^{a_{2}}=v_{2},v_{3}^{a_{1}}=v_{3},v_{3}^{a_{2}}=v_{3}v_{1}.$ Then
every $A$-orbit on $V$ has length dividing $p.$
\end{remark}

\section{Proof of theorem}

Let $(P,R,\chi)$ be a counterexample with $\left|  PR\right|  +\chi(1)$
minimum. We shall proceed in a series of steps. To simplify the notation we
set $G=PR$ .\newline

\textit{(1) $\chi$ is irreducible.}\newline

There exists an irreducible constituent $\chi_{_{1}}$ of $\chi$ which does not
contain $Z(P)$ in its kernel, that is $(\chi_{_{1}})_{_{P}}$ is faithful. Then
we have $\chi_{_{1}} = \chi$ because otherwise $\chi_{_{1}}$ contains the
regular $A$-character by induction.\newline

\textit{(2) $\chi_{_{P}}$ is homogeneous and $R_{0}=1.$ }\newline

As it is well known the irreducible characters of the extraspecial group $P$
are uniquely determined by their restriction $Z(P)$ so that $\chi_{_{P}%
}=e\theta$ for some faithful irreducible $GA$-invariant character $\theta$ of
$P$ and some positive integer $e$, since $Z(P)\leq Z(GA)$. The coprimeness
condition $(|P|,|RA_{p^{\prime}}|)=1$ enables us to extend $\theta$ in a
unique way to an irreducible character $\overline{\theta}$ of $GA_{p^{\prime}%
}$ such that $det(\overline{\theta})(x)=1$ for each $x \in RA_{p^{\prime}} $
by [\cite{Isaacs}, 8.16]. On the other hand $\theta_{1}=\theta\times1_{R_{0}}$
is an irreducible $P\times R_{0}$-character with $R_{0} \leq Ker \theta_{1}$.
We can extend $\theta_{1}$ uniquely to $\overline{\theta}_{1} \in
Irr(GA_{p^{\prime}}/R_{0})$ with $det(\overline{\theta}_{1})(x)=1$ for each $x
\in RA_{p^{\prime}}/R_{0}$. The uniqueness of this extension implies $R_{0}
\leq Ker \overline{\theta}$. Notice that $(\overline{\theta}_{1})_{_{P}%
}=\theta=\overline{\theta}_{_{P}}$ and also that the set $\{ \varphi\, : \,
\varphi\in Irr(GA_{p^{\prime}}) \,\,\text{such that} \, \varphi_{_{P}}%
=\theta\}$ is $A_{p}$-invariant, because $\theta^{a} =\theta$ for each $a \in
A_{p}$. Since $det(\overline{\theta}^{a})(x)=1 $ for each $a\in A_{p} $, the
uniqueness of $\overline{\theta}$ gives $\overline{\theta}^{a}=\overline
{\theta}$. It follows from [\cite{Isaacs}, Corollary 11.22] that
$\overline{\theta}$ is extendible to an irreducible $GA$-character, say
$\overline{\overline{\theta}}$ . Now $\overline{\overline{\theta}}_{_{G}%
}=\overline{\theta}$, $\overline{\overline{\theta}}_{_{P}}=\theta$ and $R_{0}
\leq Ker \overline{\theta}= G \cap Ker \overline{\overline{\theta}}$. If
$\overline{\overline{\theta}}(1)< \chi_{_{1}}$ or $R_{0}\ne1$, by induction
applied to the group $GA/R_{0}$ over $\overline{\overline{\theta}}$ we see
that $\overline{\overline{\theta}}_{_{A}}$ contains the regular $A$
-character. Since $\chi$ is a constituent of $\overline{\overline{\theta}%
}_{_{P}} |^{GA}$, there exists $\beta\in Irr(GA/P)$ such that $\chi
=\overline{\overline{\theta}}\cdot{\beta} $ by [\cite{Isaacs}, 6.17] and hence
$\chi_{_{A}}=\overline{\overline{\theta}}_{_{A}}\cdot\beta_{A} $. We conclude
that $\chi_{_{A}}$ contains the regular $A$-character, while $\overline
{\overline{\theta}}_{_{A}}$ does. Therefore without loss of generality we may
assume that $R_{0}=1$ as claimed.\newline

\textit{(3) Theorem follows.}\newline

Theorem 1.3 in \cite{Turull} applied to the group $PR$ over $\chi$ shows that
one of the following holds:

$(i)$ $\chi_{_{{R}}}$ contains the regular $R$-character;

$(ii)$ $p=2$, and $r$ is a Fermat prime.

By hypothesis $(d)$ we see that $(i)$ follows, that is $\chi_{_{R}}$ contains
a copy of every irreducible $R$-character. On the other hand we can regard
$Irr(R/\Phi(R))$ as a faithful $\mathbb{F}_{r}(A)$-module which is isomorphic
to $R/\Phi(R)$ and hence apply the proposition above to get a linear character
$\nu$ of $R $ such that $C_{A}(\nu)=1$. Let $V$ be a $GA$-module affording
$\chi$ and let $W$ be the homogeneous component of $V_{_{R}}$ corresponding to
$\nu$. Since the stabilizer in $A$ of $W$ is trivial, $V_{_{A}}$ contains the
regular $A$-module. Therefore $\chi_{_{A}}$ contains the regular
$A$-character. $\Box$

\end{document}